\documentstyle[11pt,leqno]{article}
\newtheorem{th}{Theorem}[section]
\newtheorem{lem}[th]{Lemma}

\title{\Large \bf {The Dolbeault operator on Hermitian spin surfaces}}
\author{{\sc B.Alexandrov} \hspace{3mm} {\sc G. Grantcharov}  \hspace{3mm} {\sc
S.Ivanov}
\thanks{The authors are supported by Contract MM 809/1998 with the
Ministry of Science and Education of Bulgaria and by Contract 238/1998 with the
University of Sofia "St. Kl. Ohridski".}}
\date{}
\begin{document}
\maketitle
\thispagestyle{empty}
\vspace{2mm}
\begin{center}
{\sc Department of Mathematics\\ University of Sofia\\ "St. Kl. Ohridski" }
\end{center}
\vspace{25mm}

\begin{abstract}
We consider the Dolbeault operator
$\sqrt 2 (\overline {\partial} + {\overline {\partial}}^*)$
of $K^{\frac{1}{2}}$ -- the square root of the canonical line bundle
which determines the spin structure of a compact Hermitian spin surface
$(M,g,J)$. We prove that all cohomology groups $H^i (M,{\mathcal O}
(K^{\frac{1}{2}}))$ vanish if the scalar curvature of $g$ is
non-negative and non-identically zero. Moreover, we estimate the first
eigenvalue of the Dolbeault operator when the conformal scalar
curvature $k$ is non-negative and when $k$ is positive. In the first case we
give a
complete list of limiting manifolds and in the second one we give non-K\"ahler
examples of limiting manifolds.
\\[15mm]
{\bf Running title:} The Dolbeault operator on Hermitian spin surfaces
\\[5mm]
\\[5mm]
${\bf MS}$ {\bf classification: } 53C15; 53C25; 53B35
\end{abstract}
\newpage

\section{Introduction}

\indent The well-known vanishing theorem of Lichnerowicz says that there are
no harmonic spinors on compact spin manifolds of non-negative
non-identically zero scalar curvature, i.e. the kernel of the Dirac
operator vanishes. When the scalar curvature is identically zero the
harmonic spinors are actually parallel. Complete
classification of the complete simply connected irreducible spin
manifolds admitting parallel spinor is given by Hitchin \cite{H} and
Wang \cite{W}. When the scalar curvature is strictly positive one may
try to find an estimate for the first eigenvalue of the Dirac
operator.  This is done by Friedrich \cite{F2}. The estimate is expressed in
terms of the scalar curvature and the limiting manifolds are characterized by
the existence of a real Killing spinor.

 It is well-known (see
\cite{H}) that in the
K\"ahler case the Dirac operator coincides with the Dolbeault operator
$\Box =\sqrt 2 (\overline {\partial} + {\overline {\partial}}^*)$
on $K^{\frac{1}{2}}$ -- the square root of the canonical line bundle
$K$ which determines the spin structure. Applying Hodge theory to the
corresponding Dolbeault complex Hitchin \cite{H} has shown that on
a compact K\"ahler spin manifold the space of harmonic spinors can
be identified with the holomorphic cohomology
$H^* (M,{\mathcal O} (K^{\frac{1}{2}}))$. So, if a compact K\"ahler spin manifold admits a
 Riemannian
metric of strictly positive scalar curvature then all cohomology groups
$H^* (M,{\mathcal O} (K^{\frac{1}{2}}))$ vanish by Lichnerowicz
theorem.

The purpose of this note is to treat problems similar to the above mentioned
in the case of the Dolbeault operator on compact Hermitian spin surfaces.

\par
Our first observation is the following.

\begin{th}\label{th1}
Let $(M,J)$ be a compact complex spin surface admitting a Hermitian
metric of non-negative non-identically zero scalar curvature. Then
$H^i (M,{\mathcal O} (K^{\frac{1}{2}}))=0$, $i=0,1,2$.
\end{th}

\noindent {\it Proof:} By arguments similar to those in Proposition~I.18 in
\cite{G4} (see also Lemma~3.3 in
\cite{ADM}) the existence of a Hermitian metric of non-negative non-identically
zero scalar curvature implies that all the plurigenera of $(M,J)$ vanish.
Hence, $H^0 (M,{\mathcal O} (K^{\frac{1}{2}}))=0$. By Serre duality
$H^2 (M,{\mathcal O} (K^{\frac{1}{2}}))=0$. By Lichnerowicz vanishing
theorem the index of the Riemannian Dirac operator $D$ vanishes. But since both
$D$ and the Dolbeault operator $\Box $ are generalized Dirac operators
on the same Clifford module, they have the same index. Thus,
$$0=ind(\Box )=dim H^0 (M,{\mathcal O} (K^{\frac{1}{2}})) -
dim H^1 (M,{\mathcal O} (K^{\frac{1}{2}})) +
dim H^2 (M,{\mathcal O} (K^{\frac{1}{2}}))$$
and therefore $H^1 (M,{\mathcal O} (K^{\frac{1}{2}}))=0$. \hfill $\Box $

In view of this result the following questions are natural:
\par
Is there  an estimate for the first
eigenvalue of the Dolbeault operator on Hermitian spin surfaces of
positive scalar curvature?
If the answer is 'Yes', describe the limiting manifolds, i.e. the
manifolds for which the estimate is attained.
\par
For K\"ahler manifolds the above questions are treated in terms of the
Dirac operator. However, according to the result of Hijazi  \cite{Hij}
the estimate of Friedrich \cite{F2}  is not
sharp on K\"ahler manifold since it admits a parallel form.
A better estimate for the
first eigenvalue of the Dirac operator on compact K\"ahler spin manifold
with strictly positive scalar curvature is found by
Kirchberg \cite{K1,K2} and the limiting manifolds are characterized by
the existense of K\"ahlerian twistor spinors.  In the 4-dimensional case
(in which we are interested in this
paper) the classification of limiting K\"ahler manifolds  is given by
Friedrich \cite{F} using detailed study of the K\"ahlerian
twistor equations. So, the
complete answer to both of the questions for compact
 K\"ahler spin surface is known.
\par
In the present paper we give complete answer to the questions
under the stronger assumption of non-negative conformal scalar
curvature. We recall that the conformal scalar curvature of a Hermitian
surface is the scalar
curvature of the corresponding Weyl connection.
In the case of strictly positive conformal scalar curvature we answer to the
first
question completely and partially to the second.

Our considerations  are based on Bochner type
calculations using the set of canonical Hermitian connections $\nabla^t, t\in
R$, described by Gauduchon \cite{G2}. Among these connections an
important role plays  the  Bismut connection. This is the unique
Hermitian connection with skew-symmetric torsion (cf. \cite{G2}) and is used
by Bismut \cite{B} to express the Weitzenb\"ock formula for the Dolbeault
operator. To treat limiting manifolds we consider twistor
equations with respect to the canonical Hermitian connections.

More precisely, we prove

\begin{th}\label{th2}
Let $(M,g,J)$ be a compact Hermitian spin surface of non-negative
conformal scalar curvature. Then the first eigenvalue $\lambda $ of
the Dolbeault operator satisfies the inequality
\begin{equation}\label{1}
\lambda ^2 \ge \frac{1}{6} {\rm inf}_M s,
\end{equation}
where $s$ is the scalar curvature of $g$. Further, the following
conditions are equivalent:

(i) There is an equality in (\ref{1}).

(ii) There exists a parallel spinor in $\Sigma _+ M$ with respect to
the Bismut connection, the conformal scalar curvature is
identically zero and the scalar curvature is constant.

(iii) $(M,g,J)$ is a $K3$-surface or a flat torus with their
hyper-K\"ahler metric or a coordinate quaternionic Hopf surface (see
\cite{Bo}) with a metric of constant scalar curvature in the conformal
class of the standard locally conformally flat metric.
\end{th}

\begin{th}\label{th3}
Let $(M,g,J)$ be a compact Hermitian spin surface of positive
conformal scalar curvature $k$. Then the first eigenvalue $\lambda $
of the Dolbeault operator satisfies the inequality
\begin{equation}\label{2}
\lambda ^2 \ge \frac{1}{2} {\rm inf}_M k.
\end{equation}

The equality in (\ref{2}) is attained iff $k$ is constant and there
exists a Hermitian twistor spinor with respect to the
Hermitian connection $\nabla ^{-3}$. In this case $(M,g,J)$ is locally
conformally K\"ahler.
\end{th}

Note that on K\"ahler surfaces the
estimate (\ref{2}) coincides with that of \cite{K2}.

In the last section we give examples of non-K\"ahler Hermitian
surfaces for which the limiting case of the inequality (\ref{2}) is
attained.

\section{Preliminaries}

Let $(M,g,J)$ be a Hermitian surface with complex structure $J$ and
compatible metric $g$. Denote by $\Omega $ the K\"ahler form, defined
by $\Omega (X,Y)=g(X,JY)$. The volume form of g is
$\omega =\frac{1}{2} \Omega \land \Omega $. It is well-known that
$d\Omega = \theta \land \Omega$, where
$\theta = \delta \Omega \circ J$ is the Lee form of $(M,g,J)$. Recall
that $(M,g,J)$ is K\"ahler iff $\theta =0$; locally conformally
K\"ahler iff $d\theta =0$; globally conformally K\"ahler iff
$\theta = df$ for a smooth function $f$ on $M$ (in this case $e^{-f}g$
is a K\"ahler metric). Let $\nabla $ be the Levi-Civita connection of
$g$ and $R$ and $s$ -- its curvature tensor and scalar curvature
respectively (for the curvature tensor we adopt the following
definition:
$R(X,Y,Z,W) = g([\nabla _X,\nabla _Y]Z - \nabla _{[X,Y]}Z,W))$.
Recall that the *-Ricci tensor $\rho ^*$ and the *-scalar curvature
$s^*$ of $M$ are defined by
$$\rho ^*(X,Y)=\sum_{i=1}^4 R(e_i,X,JY,Je_i)
=-\frac{1}{2}\sum_{i=1}^4 R(X,JY,e_i,Je_i)$$
$$s^*=\sum_{i=1}^4 \rho ^*(e_i,e_i),$$
where here and in the following $\{e_i\}$ is a local orthonormal frame
of the tangent bundle $TM$ and $\{e^i\}$ is its dual frame.

We also have (cf. \cite{V})
\begin{equation}\label{3}
s - s^* = 2\delta \theta + |\theta |^2
\end{equation}
Note that on a K\"ahler manifold the Ricci and *-Ricci tensors
coincide; in particular $s=s^*$.

Now consider the Weyl connection determined by the Hermitian structure
on $M$, i.e. the unique torsion-free connection $\nabla ^W$ such that
$\nabla ^W g = \theta \otimes g$. The conformal scalar curvature $k$
is defined to be the scalar curvature of $\nabla ^W$. Equivalently, $k$ is the
*-scalar curvature of the self-dual Weyl tensor $W_+$, multiplied by
$\frac{3}{2}$. The Weyl connection is invariant under conformal changes of
the metric since if $\widetilde{g} = e^f g$ then
$\widetilde{\theta }= \theta + df$. Hence,
\begin{equation}\label{4}
\widetilde{k} = e^{-f}k
\end{equation}
We also have
\begin{equation}\label{5}
k = \frac {3s^* - s}{2}
\end{equation}

Recall the definition of the set of canonical Hermitian connections
\cite{G2}: For a real number $t$ the connection $\nabla ^t$ is defined
by
$$g(\nabla ^t _X Y,Z) = g(\nabla _X Y,Z) - \frac{1}{2} g((\nabla
_{JX}J)Y,Z) + \frac{t}{4} (d\Omega (JX,JY,JZ) + d\Omega (JX,Y,Z))$$
or equivalently
\begin{eqnarray}\label{6}
\nabla ^t _X Y &=& \nabla _X Y - \frac{t+1}{4} \theta (Y)X +
\frac{t-1}{4} \theta (JY)JX - \frac{t}{2} \theta (JX)JY
\\ & & + \frac{t+1}{4} g(X,Y)\theta ^{\#} - \frac{t-1}{4} g(X,JY)J\theta ^{\#},
\nonumber
\end{eqnarray}
where $\theta ^{\#}$ is the vector field dual to $\theta $.

The canonical Hermitian connections form an affine line (degenerating
to a point in the K\"ahler case) determined by $\nabla ^0$ - the
projection of the Levi-Civita connection into the affine space of all
Hermitian connections, and $\nabla ^1$, which coincides with the Chern
connection. In the sequel important role will be played also by the
connections $\nabla ^{-1}$ (considered by Bismut \cite{B}) and
$\nabla ^{-3}$.

From now on we assume that $(M,g,J)$ is spin manifold. Denote by
$\Sigma M$ its spinor bundle and let
$\mu:  T^* M \otimes \Sigma M \longrightarrow \Sigma M$ be the
Clifford multiplication. Identifying $T^* M$ and $TM$ via the metric
we shall also consider $\mu $ as a map from $TM \otimes \Sigma M$ into
$\Sigma M$. We shall often write the Clifford multiplication by
juxtaposition, i.e.
$$\mu (\alpha \otimes \psi ) = \alpha \psi , \qquad
\mu (X \otimes \psi ) = X \psi , \qquad
\alpha \in T^* M, \quad X \in TM, \quad \psi \in \Sigma M.$$

Since $(M,g,J)$ is a Hermitian manifold the choice of a spin structure on
it is equivalent to the choice of a square-root $K^{\frac{1}{2}}$ of
the canonical line bundle $K$ \cite{H}. Thus for the corresponding
spinor bundle we have
\begin{equation}\label{7}
\Sigma M = \Lambda ^{0,\bullet }M \otimes K^{\frac{1}{2}},
\end{equation}
where $\Lambda ^{0,\bullet }M = \Lambda ^{0,0}M \oplus
\Lambda ^{0,1}M \oplus \Lambda ^{0,2}M$.
In particular, the spinor bundle $\Sigma M$ splits as follows:
\begin{equation}\label{8}
\Sigma M = \Sigma _0 M \oplus \Sigma _1 M \oplus \Sigma _2 M ,
\end{equation}
where
$\Sigma _r M = \Lambda ^{0,r}M \otimes K^{\frac{1}{2}}$ is the
eigensubbundle with respect to the eigenvalue $(2-2r)i$ of the
Clifford action of the K\"ahler form $\Omega $ on $\Sigma M$ (cf.
\cite{K1}).

The half-spinor bundles $\Sigma _{\pm } M$ are the eigensubbundles of
the volume form $\omega $ with respect to its eigevalues $\mp 1$ and
we have
$$\Sigma _+ M = \Sigma _0 M \oplus \Sigma _2 M, \qquad
\Sigma _- M = \Sigma _1 M.$$

Let $p_r: \Sigma M\longrightarrow \Sigma _r M$, $r=0,1,2$, be the
projections with respect to the splitting (\ref{8}). For convenience
we denote $\Sigma _{-1} M = \Sigma _3 M =0$ and $p_{-1}=p_3=0$. Recall
\cite{K1} that for $X \in TM$
\begin{equation}\label{9}
X\Omega - \Omega X = 2JX
\end{equation}
as endomorphisms of $\Sigma M$ and for $\psi \in \Sigma _r M$
\begin{equation}\label{10}
X\psi = p_{r-1} X\psi + p_{r+1} X\psi, \qquad
JX\psi = -ip_{r-1} X\psi + ip_{r+1} X\psi
\end{equation}
\begin{equation}\label{11}
p(X)\psi = p_{r+1} X\psi, \qquad \bar{p} (X)\psi = p_{r-1} X\psi ,
\end{equation}
where $p(X) = \frac{1}{2}(X-iJX)$, $\bar{p} (X) = \frac{1}{2}(X+iJX)$.

Any metric connection in the tangent bundle gives rise to metric
connection in the spinor bundle and it is easy consequence of
(\ref{6}) that
\begin{equation}\label{12}
\nabla ^t _X \psi = \nabla _X \psi + \frac{1}{8}[(t+1)X\theta +
(t-1)JXJ\theta + 2t\theta (X) + 2t\theta (JX)\Omega ]\psi
\end{equation}
for $X \in \Gamma (TM), \quad \psi \in \Gamma (\Sigma M)$.
In particular, we have
\begin{equation}\label{13}
\nabla ^t _X \psi _0 = \nabla _X \psi _0 + \frac{1}{4}X\theta \psi _0
+ \frac{1}{4}\theta (X)\psi _0 + \frac{t+1}{4}i\theta (JX)\psi _0,
\end{equation}
for $X \in \Gamma (TM), \quad \psi _0 \in \Gamma (\Sigma _0 M)$.
This follows from (\ref{12}) and
\begin{equation}\label{14}
p_0 X\alpha \psi _0 = -2\alpha (\bar{p}(X))\psi _0
\qquad  X \in \Gamma (TM), \quad  \alpha \in \Gamma (T^* M),
\quad \psi _0 \in \Gamma (\Sigma _0 M).
\end{equation}
The K\"ahler form $\Omega $ is parallel with respect to any Hermitian
connection, so the Hermitian connections preserve the splitting
(\ref{8}). Recall also \cite{ABS} that there is an antilinear bundle
map $j:\Sigma M \longrightarrow \Sigma M$ which commutes with the
Clifford multiplication by real vectors, $j^2 = -1$, $j$ is parallel
with respect to any metric connection on $M$ and preserves the
Hermitian inner product on $\Sigma M$. In particular, $j$ provides an
antilinear isomorphism between $\Sigma _0 M$ and $\Sigma _2 M$.

The Dirac operator of the Levi-Civita connection
$D:\Gamma (\Sigma M) \longrightarrow \Gamma (\Sigma M)$ is defined by
$D=\mu \circ \nabla $ or equivalently by
$D\psi = \sum_{i=1}^4 e^i \nabla _{e_i}\psi $. The Dirac operators
$D^t$ of the connections $\nabla ^t$ are defined in similar way by
replacing $\nabla $ with $\nabla ^t$.

The identification (\ref{7}) shows that the Dolbeault operator
$\Box =\sqrt 2 (\overline {\partial} + {\overline {\partial}}^*)$
of $K^{\frac{1}{2}}$ also acts on sections of the spinor bundle. As
shown in \cite{G2}
\begin{equation}\label{15}
D^t = D - \frac{3t}{4}\theta - \frac{2t-1}{4}J\theta \circ \Omega
\end{equation}
\begin{equation}\label{16}
\Box = D + \frac{1}{4}\theta + \frac{1}{4}J\theta \circ \Omega
\end{equation}
Let $\pi_r:T^* M\otimes \Sigma _r M \longrightarrow Ker\mu
\vert_{T^* M \otimes \Sigma _r M}$, $r=0,1,2$, be the orthogonal
projections.  It is easy to see that
\begin{eqnarray}\label{17}
\pi _r(\alpha \otimes \psi ) & &
\\ = \alpha \otimes \psi &+& \frac{1}{2(r+1)} \sum_{i=1}^4 e^i \otimes
\bar{p}(e^i)p(\alpha )\psi + \frac{1}{4-2(r-1)}\sum_{i=1}^4 e^i \otimes
p(e^i)\bar{p}(\alpha )\psi , \nonumber
\end{eqnarray}
where $p(\alpha ) = \frac{1}{2}(\alpha -iJ\alpha )$,
$\bar{p} (\alpha ) = \frac{1}{2}(\alpha +iJ\alpha )$ and $J\alpha $ is
the dual form of $J\alpha ^{\#}$, $J\alpha = - \alpha \circ J$.

The twistor operators of the Hermitian connection $\nabla ^t$ are the
differential operators \\
$P^t_r:~\Gamma (\Sigma _r M)~\longrightarrow ~
\Gamma (T^* M \otimes \Sigma _r M)$ defined by
$P^t_r = \pi _r\circ \nabla ^t$. By analogy with the K\"ahler twistor
spinors of \cite{K2} we shall call the spinors in the kernel of $P^t_r$
{\it Hermitian twistor spinors with respect to} $\nabla ^t$. We are
particularly interested in $P^t_0$. It follows from (\ref{17}) that
\begin{equation}\label{18}
P^t_0 \psi _0 = \nabla ^t \psi _0 + \frac{1}{2}
\sum_{i=1}^4 e^i \otimes p_0 e_i D^t \psi _0, \qquad
\psi _0 \in \Gamma (\Sigma _0 M)
\end{equation}
and also
\begin{equation}\label{19}
P^t_0 \psi _0 = \sum_{i=1}^4 e^i \otimes \nabla ^t_{p(e_i)}\psi _0 =
(\nabla ^t)^{1,0} \psi _0, \qquad \psi _0 \in \Gamma (\Sigma _0 M).
\end{equation}
In the 4-dimensional case the Weitzenb\"ock formula for the Dolbeault
operator (\cite{B}, Theorem~2.3) reads as follows:
\begin{equation}\label{20}
\Box ^2 = \Delta ^{-1} + \frac{s}{4} + \frac{1}{4}\delta \theta .\omega
- \frac{|\theta |^2}{8},
\end{equation}
where $\Delta ^t = (\nabla ^t)^* \nabla ^t$ is the spinor Laplacian of
the connection $\nabla ^t$.

In the following we shall denote by $<.,.>$ and $|.|$ pointwise inner
products and norms and by $(.,.)$ and $\| .\|$ -- the global ones
respectively.

\section{Proof of Theorem~\ref{th2} and Theorem~\ref{th3}}
Throughout this section $(M,g,J)$ is a compact Hermitian spin surface.

\begin{lem}\label{lem1}
Let $\lambda \not = 0$ be an eigenvalue of the Dolbeault operator
$\Box $. Then there exists an eigenspinor
$\psi _0 \in \Gamma (\Sigma _0 M)$ for the eigenvalue $\lambda ^2$ of
$\Box ^2$.
\end{lem}

\noindent {\it Proof:} It is clear that
$\Box (\Gamma (\Sigma _r M)) \subset \Gamma (\Sigma _{r-1} M) \oplus
\Gamma (\Sigma _{r+1} M)$, and
$\Box ^2 (\Gamma (\Sigma _r M)) \subset \Gamma (\Sigma _r M)$,
$r=0,1,2$. Hence, if $\psi \not = 0$ is an eigenspinor of $\Box $ for
$\lambda $ and $\psi = \psi _0 + \psi _1 + \psi _2$ with respect to
the decomposition (\ref{8}) then
\begin{equation}\label{21}
\psi = (\psi _0 + \frac{1}{\lambda }\Box \psi _0) +
(\psi _2 + \frac{1}{\lambda }\Box \psi _2)
\end{equation}
(i.e. $\psi _1 = \frac{1}{\lambda }\Box \psi _0 +
\frac{1}{\lambda }\Box \psi _2$) and both
$\psi _0 + \frac{1}{\lambda }\Box \psi _0$ and
$\psi _2 + \frac{1}{\lambda }\Box \psi _2$ are eigenspinors for
$\Box$ with respect to $\lambda $. Moreover,
$\psi _0$, $\Box \psi _0$, $\psi _2$, $\Box \psi _2$
are eigenspinors for $\Box ^2$ with respect to $\lambda ^2$.
It follows by
(\ref{16}) that $\Box \circ j = j \circ \Box $ and therefore
$j\psi _2 \in \Gamma (\Sigma _0 M)$ is also an eigenspinor for
$\Box ^2$ with respect to $\lambda ^2$.
Since by (\ref{21}) $\psi _0 \not = 0$ or $\psi _2 \not = 0$,
Lemma~\ref{lem1} is proved. \hfill $\Box $

\vspace{3mm}
\noindent {\it Proof of Theorem~\ref{th2}:} When
${\rm inf}_M s \le 0$ the inequality (\ref{1}) is trivially satisfied.
So we can assume that ${\rm inf}_M s > 0$. Hence it follows from
Theorem~\ref{th1} that $\lambda \not = 0$ and by Lemma~\ref{lem1} we
have an eigenspinor $\psi _0 \in \Gamma (\Sigma _0 M)$ of $\Box ^2$
with respect to $\lambda ^2$. Thus (\ref{20}) yields
$$\lambda ^2 \| \psi _0 \| ^2=(\Box ^2 \psi _0,\psi _0) =
(\Delta ^{-1} \psi _0,\psi _0) + \frac{1}{4}
((s - \delta \theta - \frac{|\theta |^2}{2}) \psi _0,\psi _0).$$
By (\ref{3}) $s - \delta \theta - \frac{|\theta |^2}{2} =
\frac{s+s^*}{2}$ and since $k \ge 0$, i.e. $s^* \ge \frac{s}{3}$, we
obtain that $\lambda ^2 \ge \frac{1}{6} {\rm inf}_M s$.

\noindent Now we proceed with the proof of the second part of the theorem. Note
that the Bismut connection restricted to sections of $\Sigma _+ M$
coincides with the connection considered in \cite{M} and called there
"the Weyl connection". Hence, as proved in \cite{M}, a parallel spinor
in $\Sigma _+ M$ with respect to the Bismut connection gives rise to a
hyper-Hermitian structure on $M$ and thus $M$ is conformally
equivalent to one of the manifolds in (iii) (cf. also \cite{Bo}). In
particular, it follows that $M$ is anti-self-dual, which is equivalent
to $k=0$ and $d\theta = 0$ (cf. for example \cite{ADM}). Conversely,
any manifold conformally equivalent to those listed in (iii) admits a
parallel spinor in $\Sigma _+ M$ (and hence in $\Sigma _0 M$) with
respect to the Bismut connection in the spin structure given by the
trivial square-root of the canonical line bundle. This proves the
equivalence of (ii) and (iii) since the hyper-K\"ahler metrics on a
$K3$-surface or a torus are the only metrics of constant scalar
curvature in their conformal classes.

\noindent Now we prove the equivalence of (i) and (ii).

\noindent When $s>0$ the calculations in the first part of the proof show that
there is an equality in (\ref{1}) iff $s=const$, $k=0$ and
$\nabla ^{-1} \psi _0 = 0$. But as shown above,
$\nabla ^{-1} \psi _0 = 0$ implies $k=0$ and thus (i) and (ii) are
equivalent when $s>0$.

\noindent When ${\rm inf}_M s = 0$ Theorem~\ref{th1} shows that the equality in
(\ref{1}) is possible only if $s=0$. Hence, by $k \ge 0$ we have
$s^* \ge 0$ and integrating (\ref{3}) we obtain that $\theta =0$, i.e.
$(M,g,J)$ is K\"ahler. Therefore (\ref{20}) (which in the K\"ahler
case coincides with the usual Lichnerowicz formula and the Bismut
connection coincides with Levi-Civita connection) yields that
$\Box ^2 \psi = 0$ iff $\nabla ^{-1} \psi = 0$. Hence it remains to
show that if there exists a parallel spinor, then there exists
a parallel spinor in $\Sigma _+ M$. But $(M,g,J)$ is anti-self-dual
(since $k=s=0$ and $d\theta = 0$) and thus its signature
$\sigma (M) \le 0$. Since $\sigma (M) = -8ind(D)=-8ind(\Box )$ it
follows that $ind(\Box ) \ge 0$, i.e.
$dim Ker(\Box \vert _{\Sigma _+ M}) \ge dim Ker(\Box \vert _{\Sigma _- M})$.

\noindent Thus the equivalence of (i) and (ii) is proved. \hfill $\Box $

We start the proof of Theorem~\ref{th3} by two Lemmas.

\begin{lem}\label{lem2}
Let $\psi _0 \in \Gamma (\Sigma _0 M)$. Then
\begin{eqnarray}\label{22}
\| P^t_0 \psi _0 \| ^2 &=& \frac{1}{2}(\Box ^2 \psi _0 ,\psi _0 ) -
\frac{1}{4}(s \psi _0 ,\psi _0 ) -
\frac{t}{4}(\delta \theta \psi _0 ,\psi _0 ) \\
& & + \frac{t^2 - 2t -3}{32}(|\theta | ^2 \psi _0 ,\psi _0 ) -
\frac{t+3}{4}\Re (\theta \Box  \psi _0 ,\psi _0 ). \nonumber
\end{eqnarray}
\end{lem}

\noindent {\it Proof:} It follows by (\ref{18}) that
$$| P^t_0 \psi _0 | ^2 = | \nabla ^t \psi _0 | ^2 -
| D^t \psi _0 | ^2  -
\frac{1}{4}\sum _{j=1}^{4}
<e_j \bar{p} (e_j) D^t \psi _0 , D^t \psi _0 >.$$
But $\sum _{j=1}^{4}e_j \bar{p} (e_j) = \frac{1}{2} (-4 +
i\sum _{j=1}^{4}e_j Je_j ) = -2 - i\Omega $. Thus
$\sum _{j=1}^{4}e_j \bar{p} (e_j) D^t \psi _0 = -2D^t \psi _0$ since
$D^t \psi _0 \in \Gamma (\Sigma _1 M)$. Hence,
\begin{equation}\label{23}
| P^t_0 \psi _0 | ^2 = | \nabla ^t \psi _0 | ^2 -
\frac{1}{2}| D^t \psi _0 | ^2.
\end{equation}
By (\ref{15}) and (\ref{16}) we obtain
\begin{equation}\label{24}
D^t \psi _0 = \Box \psi _0 + \frac{t-1}{4} \theta \psi _0
\end{equation}
and hence
\begin{equation}\label{25}
|D^t \psi _0 | ^2 = |\Box \psi _0 | ^2 +
\frac{(t-1)^2}{16} <|\theta | ^2 \psi _0,\psi _0> -
\frac{t-1}{2}\Re <\theta \Box  \psi _0 ,\psi _0 >.
\end{equation}
By (\ref{13}) we obtain
\begin{equation}\label{26}
\nabla ^t _X \psi _0 = \nabla ^{-1}_X \psi _0
+ \frac{t+1}{4}i\theta (JX)\psi _0,
\end{equation}
Hence,
\begin{equation}\label{27}
|\nabla ^t \psi _0 | ^2 = |\nabla ^{-1} \psi _0 | ^2 +
\frac{(t+1)^2}{16} <|\theta | ^2 \psi _0,\psi _0> +
\frac{t+1}{2}\Re <i\nabla ^{-1}_{J\theta ^{\#}} \psi _0 ,\psi _0 >.
\end{equation}
It is easily seen that $D^t \circ \alpha + \alpha \circ D^t =
\sum _{j=1}^{4} e^j \nabla ^t _{e_j} \alpha -2\nabla ^t _{\alpha ^{\#}}$
for arbitrary 1-form $\alpha $. Hence
\begin{eqnarray}
<\nabla ^t_{J\theta ^{\#}} \psi _0 &,& \psi _0 > \nonumber \\
&=& \frac{1}{2}\sum _{j=1}^{4}
<e^j \nabla ^t _{e_j} (J\theta ) \psi _0 ,\psi _0 > -
\frac{1}{2}<D ^t \circ J\theta \psi _0 ,\psi _0 > -
\frac{1}{2}<J\theta \circ D ^t \psi _0 ,\psi _0 > \nonumber \\
&=& \frac{i}{2}\sum _{j=1}^{4}
<e^j \nabla ^t _{e_j} \theta \psi _0 ,\psi _0 > -
\frac{i}{2}<D ^t \circ \theta \psi _0 ,\psi _0 > +
\frac{i}{2}<\theta \circ D ^t \psi _0 ,\psi _0 > \nonumber \\
&=& i<\theta \circ D ^t \psi _0 ,\psi _0 > +
i<\nabla ^t_{\theta ^{\#}} \psi _0 ,\psi _0 >. \nonumber
\end{eqnarray}
Substituting this
equation in (\ref{27}) and using (\ref{24}) and the fact that
$$\Re <\nabla ^t_{\theta ^{\#}} \psi _0 ,\psi _0 > =
\frac{1}{2}(<\nabla ^t_{\theta ^{\#}} \psi _0 ,\psi _0 > +
<\psi _0 ,\nabla ^t_{\theta ^{\#}} \psi _0 >) =
\frac{1}{2}\theta ^{\#} (|\psi _0 | ^2) =
\frac{1}{2}g(\theta ,d(|\psi _0 | ^2)),$$
we obtain
\begin{eqnarray}\label{28}
|\nabla ^t \psi _0 | ^2 &=& |\nabla ^{-1} \psi _0 | ^2 +
\frac{t^2 - 2t - 3}{16} <|\theta | ^2 \psi _0,\psi _0> -
\frac{t+1}{4}g(\theta ,d(|\psi _0 | ^2)) \\
& & - \frac{t+1}{2}\Re <\theta \Box \psi _0 ,\psi _0 >. \nonumber
\end{eqnarray}
Now (\ref{22}) is obtained by substituting (\ref{25}) and (\ref{28})
in (\ref{23}), integrating and using (\ref{20}) to express
$\| \nabla ^{-1} \psi _0 \| ^2 = (\Delta ^{-1} \psi _0 ,\psi _0 ) =
(\Box ^2 \psi _0 ,\psi _0) - \frac{1}{4}
((s - \delta \theta  - \frac{|\theta |^2}{2})\psi _0,\psi _0)$. \hfill
$\Box $

\begin{lem}\label{lem3}
Let $P^t_0 \psi _0 = 0$, where $\psi _0 \in \Gamma (\Sigma _0 M)$ is
non-identically zero and $t \not = 1$. Then $(M,g,J)$ is locally
conformally K\"ahler.
\end{lem}

\noindent {\it Proof:} By (\ref{19}) we have
$\nabla ^t _Z \psi _0 = 0$ for $Z \in T^{1,0} M$. Hence, the curvature
$R^t(Z_1 , Z_2 )\psi _0 = 0$ for $Z_1 ,Z_2 \in T^{1,0} M$. It follows
by Lemma~\ref{lem2} that the principal symbol of $(P^t_0)^* P^t_0$ is
multiple of identity and thus by theorem of Aronszjan \cite{A}
$\psi _0 \not = 0$ on dense open subset of $M$. Since $\Sigma _0 M$ is
line bundle, it follows that the curvature form $R^t$ is $(1,1)$-form
on this dense open subset and hence on the whole $M$. By (\ref{13}) we
obtain $R^t = R^1 -\frac{t-1}{4}id(J\theta )$. But $R^1$ is the
curvature of the Chern connection and hence is also a $(1,1)$-form.
Thus $dJ\theta $ is $(1,1)$-form and therefore $d\theta $ is
$(1,1)$-form. Since $d\theta \bot \Omega$ it follows that $d\theta $
is anti-self-dual and the compactness of $M$ implies that
$d\theta = 0$ , i.e. $(M,g,J)$ is locally conformally K\"ahler. \hfill
$\Box $

\vspace{3mm}
\noindent {\it Proof of Theorem~\ref{th3}:} By theorem of Gauduchon
\cite{G3} there always exists a metric $\widetilde{g}$ in the conformal
class of $g$ whose Lee form is co-closed. By (\ref{4}) $\widetilde{k} >0$
and by (\ref{3}) and (\ref{5}) it follows that $\widetilde{s} >0$. Hence,
Theorem~\ref{th1} tels us that $Ker(\Box )=\lbrace 0 \rbrace $, i.e.
$\lambda \not = 0$. Thus, by Lemma~\ref{lem1} there exists an
eigenspinor $\psi _0 \in \Gamma (\Sigma _0 M)$ of $\Box ^2$ with
respect to $\lambda ^2$. Applying Lemma~\ref{lem2} with $t=-3$ we
obtain
$$(\Box ^2 \psi _0 ,\psi _0 ) = 2\| P^{-3}_0 \psi _0 \| ^2 +
\frac{1}{2}((s-3\delta \theta - \frac{3}{2}|\theta | ^2) \psi _0 ,
\psi _0 ).$$ But it follows by (\ref{3}) and (\ref{5}) that
$s-3\delta \theta - \frac{3}{2}|\theta | ^2 = k$. Thus
$$\lambda ^2 \| \psi _0 \| ^2 = 2\| P^{-3}_0 \psi _0 \| ^2 +
\frac{1}{2}(k\psi _0 ,\psi _0 ).$$ Hence,
$\lambda ^2 \ge \frac{1}{2} {\rm inf}_M k$ and the equality is
attained iff $P^{-3}_0 \psi _0 =0$ and $k=const$.

\noindent The last statement of the theorem follows from Lemma~\ref{lem3}.
\hfill $\Box $

{\bf Remark:} In Lemma~\ref{lem3} we characterize the compact Hermitian spin
surfaces
admitting Hermitian twistor spinor with respect to any of the canonical
Hermitian connections except the Chern connection.  In the same direction, one
can see that under the additional assumptions of $b_1$ odd and positive
fundamental constant $C(M,g)$ (see the definition in \cite{G1}) a compact
Hermitian spin surface admitting Hermitian twistor spinor with respect to the
Chern connection is biholomorphic to a primary Hopf surface.  To prove this
notice that the positivity of $C(M,g)$ implies by Gauduchon's plurigenera
theorem \cite{G1} that all the plurigenera of $(M,J)$ vanish and since $b_1$ is
odd $(M,J)$ must be of type $VII_0$.  A Hermitian twistor spinor with respect
to the Chern connection is an antiholomorphic section $\psi $ of
$K^{\frac{1}{2}}$.  Equivalently, $j\psi $ is holomorphic section of
$K^{-\frac{1}{2}}$ and hence $H^0 (M,{\mathcal O} (K^{-1}))\not =0$.  Now the
statement follows by the fact that $M$ is spin and arguments similar to those
in Corollary~3.12 in \cite{P}.

\section{Examples}

In this section we give non-K\"ahler examples of Hermitian surfaces
for which the limiting case in Theorem~\ref{th3} is attained.

Recall (cf. \cite{BFGK}) that under a conformal change of the metric
$\widetilde{g} = e^f g$ of a spin manifold $(M,g)$ there is an
identification $\quad \widetilde { } \quad $ of the spinor bundle $\Sigma M$ of
$(M,g)$
and the spinor bundle $\widetilde{\Sigma} M$ of $(M,\widetilde{g})$ such that
\begin{equation}\label{29}
\widetilde{X\psi } = \widetilde{X} \widetilde{\psi }
\end{equation}
\begin{equation}\label{30}
\widetilde{\alpha \psi } = \widetilde{\alpha } \widetilde{\psi }
\end{equation}
\begin{equation}\label{31}
\widetilde{\nabla } _X \widetilde{\psi } = \widetilde{\nabla _X \psi } -
\frac{1}{4} Xdf\widetilde{\psi } - \frac{1}{4} Xf.\widetilde{\psi },
\end{equation}
where $\psi \in \Gamma (\Sigma M)$,
$\widetilde{\psi } \in \Gamma (\widetilde{\Sigma } M)$ is the spinor
corresponding to $\psi $, $X \in TM$, $\widetilde{X} = e^{-\frac{f}{2}} X$,
$\alpha \in T^* M$, $\widetilde{\alpha } = e^{\frac{f}{2}} \alpha $, and
$\nabla $ and $\widetilde{\nabla}$ are the Levi-Civita connections of $g$
and $\widetilde{g}$ respectively.

Now, if $(M,g,J)$ is a Hermitian surface, the K\"ahler form of
$(M,\widetilde{g},J)$ is $\widetilde{\Omega } = e^f \Omega $ and it follows
from (\ref{30}) that
\begin{equation}\label{32}
\widetilde{\Omega \psi } = \widetilde{\Omega } \widetilde{\psi }.
\end{equation}
Thus the eigensubbundles $\Sigma _r M$ of $\Omega $ correspond
to the eigensubbundles $\widetilde{\Sigma } _r M$ of $\widetilde{\Omega }$. The
Lee form of $\widetilde{g}$ is $\theta + df$ and hence by (\ref{13}),
(\ref{29})-(\ref{32}) we obtain
\begin{equation}\label{33}
\widetilde{\nabla } ^t_X \widetilde{\psi _0} = \widetilde{\nabla ^t_X \psi _0} -
\frac{t+1}{4} i(Jdf)(X)\widetilde{\psi _0}, \qquad
\psi _0 \in \Gamma (\Sigma _0 M)
\end{equation}

\begin{lem}\label{lem4}
Let $(M,g,J)$ be a Hermitian spin surface and $\widetilde{g} =e^f g$.
If $\psi _0 \in \Gamma (\Sigma _0 M)$
is a Hermitian twistor spinor with respect to $\nabla ^t$ then
$e^{\frac{t+1}{4} f} \widetilde{\psi _0}$ is a Hermitian twistor spinor with
respect to $\widetilde{\nabla} ^t$.
\end{lem}

\noindent {\it Proof:} Easy consequence of (\ref{19}) and (\ref{33}).
\hfill $\Box $

Now let $M$ be one of the manifolds $S^2 \times S^2$ and
$T^2 \times S^2$, where $S^2$ is the 2-dimensional sphere and $T^2$
is 2-dimensional flat torus. With their standard product metrics these
manifolds are limiting K\"ahler manifolds, i.e. they are K\"ahler and
(\ref{2}) is equality for them (cf. \cite{K2} or \cite{F}). The
limiting K\"ahler manifolds are characterized by constant positive
scalar curvature and existence of K\"ahlerian twistor spinor, i.e. an
antiholomorphic section of $K^{\frac{1}{2}}$ - the square root of the
canonical line bundle which determines the spin structure (cf.
\cite{K3} or Theorem~\ref{th3} above).

We can change conformally the metric of the first factor of $M$ so to
obtain a metric $g$ on $M$ of positive non-constant scalar curvature
$s$. This metric will be K\"ahler with respect to the same complex
structure $J$ and $M$ will still admit K\"ahlerian twistor spinor.
Now we change the metric conformally by $\widetilde{g} = sg$. Since in the
K\"ahler case the connections $\nabla ^t$ coincide with the
Levi-Civita connection, it follows from Lemma~\ref{lem4} that on
$(M,\widetilde{g} ,J)$ there exists a Hermitian twistor spinor with
respect to $\nabla ^t$ (and in particular $\nabla ^{-3}$) and by
(\ref{4}) the conformal scalar curvature is constant:
$\widetilde{k} = 1$. Thus by Theorem~\ref{th3} the inequality (\ref{2}) turns
into equality for the non-K\"ahler Hermitian surface $(M,\widetilde{g} ,J)$.

\bigskip {\bf Authors' permanent address}\\[2mm] Bogdan Alexandrov, Gueo
Grantcharov, Stefan Ivanov\\ University of Sofia, Faculty of
Mathematics and Informatics, Department of Geometry, \\ 5 James
Bourchier blvd, 1126 Sofia, BULGARIA.\\ E-mail: B.Alexandrov:
alexandrovbt@fmi.uni-sofia.bg \par
 G.Grantcharov: geogran@fmi.uni-sofia.bg \par
S.Ivanov: ivanovsp@fmi.uni-sofia.bg \\
Present address of Gueo Grantcharov:\\
Department of Mathematics, University of Californiya at Riverside, Riverside,
CA 92521.
\end{document}